    \documentclass [12pt]
{amsart}
    \usepackage{amsfonts}
    \usepackage{amssymb}
    \usepackage{amsbsy}
\usepackage{amsthm}
\usepackage{amsmath}

\makeatletter

\newcommand{\dom}{\operatorname{dom}}

\newcommand{\on}{{\upharpoonright}}

\renewcommand{\cases}[1]{
\left \{\,\vcenter
 {\normalbaselines \m@th \ialign {$##\hfil $&\quad ##\hfil
 \crcr #1\crcr }}\right .}
\makeatother

\swapnumbers

\def\xx#1 {%
\newtheorem{#1}[thm]{#1}}
\let\yy\xx

\theoremstyle{plain}
\xx Theorem
\xx Corollary
\xx Lemma
\xx {Main Lemma}
\xx {Independence Lemma}
\xx Claim
\xx Conclusion
\yy Fact
\xx {Canonization Theorem}

\theoremstyle{definition}
\yy Definition
\yy Notation
\yy Setup
\yy {Fact and Definition}

\theoremstyle{remark}

\yy{Case 1}
\yy{Case 2}
\yy Remark

\yy Acknowledgments
\yy Note
\yy {Overview of the construction}

\newenvironment{Proof}{\begin{proof}}{\end{proof}}

\def\itm#1 {\item[(#1)]}

\newcommand{\RR}{{\mathbf R}}

\begin{document}

\title{Order polynomially complete lattices must be LARGE}

\author{Martin Goldstern}
\thanks{The first author is supported by an Erwin  Schr\"odinger 
fellowship from  the Austrian Science Foundation}
\address{Algebra und Diskrete Mathematik\\
Technische Universit\"at Wien\\
Wiedner Hauptstra\ss e 8-10/118.2\\
A-1040 Wien, Austria}
\curraddr{Mathematik WE 2\\
Freie Universit\"at\\
Arnimallee 3\\
D-14195 Berlin, Germany}
\email{Martin.Goldstern@tuwien.ac.at}

\author{Saharon Shelah}
\thanks{The second author is supported by the
   German-Israeli Foundation for Scientific Research \& Development
   Grant No. G-294.081.06/93.   Publication number 633.}
\address{Department of Mathematics\\
Hebrew University of Jerusalem\\
Givat Ram\\
91904 Jerusalem, Israel}
\email{shelah@math.huji.ac.il}

\subjclass{Primary 06A07; 
%       Combinatorics of partially ordered sets
%
secondary
  08A40, 
%       Operations, polynomials, primal algebras
%
  06B99,
%       None of the above but in this section (06Bxx=Lattices)
% 
03E55} 
%       Large cardinals

%%%  03E05   
%%%%    Combinatorial set theory, See also {04A20}
%%%  03E35  
%%%%    Consistency and independence results


\date{July 4, 1997}

\begin{abstract}
If $L$ is an o.p.c.\ (order polynomially complete) lattice,
then the cardinality of $L$ is a strongly inaccessible
cardinal.  In particular, the existence of o.p.c.\ lattices is
not provable in ZFC, not even from ZFC+GCH.
\end{abstract}

\maketitle

% 
% \makeatletter
% \let\old@item\@item
% \def\itemizeitem{\@ifnextchar[{\itemize@item}{\itemize@item[\@itemlabel]}}
% \def\itemize@item[#1]{\old@item[{#1}]\if@noparitem\else
%                         \def\@currentlabel{#1}\fi}
% \makeatother

\section*{Introduction}

Let $(L_1, \le) $ and $(L_2, \le)$ be partial orders.  We call a
map $f:L_1 \to L_2$ ``monotone'' if $x \le y$ implies $f(x) \le f(y)$.

Notice that if $(L, \wedge, \vee, \le)$ is a lattice, then every
polynomial function (i.e., function induced by a 
lattice-theoretic polynomial) is monotone.

% Let $(L, {\logand}, {\vee})$ be a lattice.  The  {\bf polynomials}
% over $L$ are the well-formed expressions involving constants
% from $L$ and/or   formal variables $\x_1$, $\x_2$, \dots\
%  together with parentheses and
% the formal symbols $\logand$, $\vee$.  The set of
% $n$-ary {\bf  polynomial functions}
% is the smallest set of functions from $L^n$ to $L$ which
% contains all the constant functions  as well as  the projection functions
% and is closed under (pointwise) meet and join. Every polynomial in
% which at most the variables $\x_1$, \dots, $\x_n$ occur 
% defines in an obvious way an $n$-ary polynomial function. Clearly every
% polynomial function $f$  is monotone. 

We call a lattice $L$ {\bf $n$-order polynomially complete (o.p.c.)}  if
\begin{itemize}
\itm $*$  every monotone $n$-ary function is a polynomial function 
\end{itemize}
and we say that $L$ is order polynomially complete if $L$ is $n$-order
polynomially complete for every $n$.

\medskip

The question whether there can be an infinite o.p.c.\  lattice has been the
subject of several papers. Kaiser and Sauer \cite{KS} remarked that
such a lattice cannot be countable, and 
 Haviar and
Plo\v s\v cica showed in \cite{HP97} that such a lattice would have to
be at least of size $\beth_\omega$.   (Here, $\beth_0 = {\aleph_0}$,
$\beth_{n+1} = 2^{\beth_n}$, $\beth_\omega = \sup\{\beth_n: n \in
\omega \}$. 

We show here that the size of an infinite o.p.c.\  lattice (if one exists
at all) must be a strongly inaccessible cardinal.

  In particular, the
existence of such a lattice cannot be derived from the ``usual'' axioms
of mathematics, as codified in the Zermelo-Fraenkel axioms for set
theory.   Moreover, also certain  additional assumptions
such as the (generalized) continuum hypothesis are not sufficient to
prove the existence of an o.p.c.\  lattice, or in other words, the theory
ZFC+GCH+``there is no infinite o.p.c.\  lattice'' is 
consistent\footnote{Pedants are invited to insert the 
necessary disclaimer themselves.}.
In fact, \addtocounter{footnote}{-1}%
the well-known consistent\footnotemark\ 
  theory ZFC+``there is no inaccessible cardinal'',
a natural extension of ZFC,  proves that there is no o.p.c.\  lattice.

We still do not know whether the existence of an o.p.c.\  lattice can be
refuted in ZFC alone.

% 
% 
%   Unfortunately (for
% algebraists reading this paper) or fortunately (for logicians) there
% is very little algebra in this paper.   The main tools used are
% several partition theorems, and some infinitary ``counting''
% arguments. 
% The bulk of this paper is devoted to showing that the cardinality of
% an opc lattice cannot be a singular cardinal 

In \cite{GoSh:554} we showed that if we change the  original 
question by relaxing ``lattice'' to ``partial order'', and
``polynomial'' to ``definable'', then we get (consistently) a positive
answer, already for a partial order of size $\aleph_1 $.

\noindent{\bf Acknowledgement.}  We are grateful to Lutz Heindorf
for his thoughtful comments on an earlier version of this paper, and
for alerting us to \cite{HP97}.

\section{Preliminaries}

We define here some of the notation that we will use, and we quote
several well-known theorems and corollaries 
from the calculus of partitions.

\begin{Definition}\label{rdef}
We fix a set ${\RR}   =
 \{ {\boldsymbol <}, {\boldsymbol >} , {\boldsymbol =}, {\boldsymbol
 \|} \} $ of 4 
symbols.
% , partially ordered in the natural way.  (That is, $\boldsymbol>$ is the
%greatest element, $\boldsymbol<$ the smallest, and $\boldsymbol=$ and
% $\boldsymbol\|$ are incomparable.)
 For any p.o.\ $(L, \le)$ we define 
$R: L \times L \to {\RR}$ in the obvious way: 
$ R(x,y) = {\boldsymbol<} $ iff $x<y$, etc. 
\end{Definition}

\begin{Definition}\label{defis}
 Let $(L , \le)$ be a partial order.
\begin{enumerate}
\item    We say that a set $A
\subseteq L$ is ``co-well-ordered'' iff $(A, \ge)$ is a well-ordered
set.   
\item    We call a set $A \subseteq L$
``uniform'' iff $A$ is either an antichain, or a well-ordered chain, or a
co-well-ordered chain.   
\item 
If $A$ is  well-ordered we say that the type of $A$ is 
``$\boldsymbol<$''. Similarly
we define the types ``$\boldsymbol >$'' and ``$\boldsymbol \|$''. 
[This notation is ambiguous
if $A$ is finite.  However, we are mainly interested in (large)
infinite sets anyway, so this ambiguity will not cause any problems.]
\item 
Let  $A_i \subseteq L$ for $i \in I $.
We call $(A_i:i\in I)$
``canonical'' if the following conditions hold: 
\begin{enumerate}
\item  each $A_i$ is uniform, all $A_i$ are of the same type, 
\item there is a function $F:I \times I \to  \RR$
such that: 
$$ \forall i,j\in I: i\not=j\Rightarrow
\forall a\in A_i\,\, \forall b\in A_j: R(a,b)=F(i,j).$$
Loosely speaking, this says that whenever $i\not=j$, then either $A_i$
lies ``completely above'' $A_j$, or conversely, or $A_i$ is
``completely incomparable'' with $A_j$. 
\end{enumerate}
\end{enumerate}
\end{Definition}

\begin{Definition}
For any set $A$, $[A]^2$ is the set of unordered pairs from $A$: 
$$ [A]^2 = \{ \{x,y\}: x,y\in A, x\not=y\}$$
When we consider a set $A$ together with an 1-1 enumeration 
$A = \{ a_i: i < \kappa \} $, we usually identify the unordered 
pair  $\{ a_i, a_j \} \in [A]^2$ with the ordered pair 
$(a_i, a_j)$ whenever $i<j$. 
\end{Definition}

\begin{Definition}
Let $\kappa $, $\lambda $ and $c$ be cardinals.
  The ``partition symbol''
$$ \lambda \to (\kappa)^2_c$$
means:  Whenever $F:[L]^2 \to C$, where $|L|= \lambda$, 
$|C|=c$, then there is an ``$F$-homogeneous set''
 $K \subseteq L$ of cardinality $ \kappa$, i.e., a set $K$
such that  $F \on [K]^2$ is constant.

(In other words,  whenever the edges of the complete
graph on $\lambda$ many vertices are
colored with $c$ colors, then there is a complete subgraph with
$\kappa$ many vertices all of whose edges have the same color.)
\end{Definition}

\begin{Definition}\label{munudef}
Let $(L,\le)$ be a partial order.   We will try to get some
information on the structure of $L$ by considering certain  ``cardinal
characteristics'' $\mu(L)$ and $\nu(L)$, which are defined as
follows: 
\begin{enumerate}
\item 
        We let $\mu(L)$ be the smallest cardinal $\mu$ such that there
        is no uniform set $A \subseteq L$ of cardinality $\mu$.  In
        other words, $\kappa < \mu(L)$ iff there is a uniform subset
        $A \subseteq L$ of size $\kappa$.

\item We let $\mu_n(L)= \mu(L^n)$ for $n>0$.
\item
We let $\nu(L)$ be the smallest cardinal $\nu$ such that there is no
family $(f_i: i<\nu)$ of $\nu$ many pairwise incomparable monotone 
functions from $L$ to $L$.     (Functions are ordered pointwise.)  
\item 
$\nu_n(L)$ is the smallest cardinal $\nu$ such that there are no
pairwise incomparable monotone functions $(f_i:i<\nu)$ from $L^n$ to $L$.
\item More generally, $\nu(L_1, L_2)$ is the
 the smallest cardinal $\nu$ such that there are no
pairwise incomparable monotone functions $(f_i:i<\nu)$ from $L_1 $ to $L_2$.
\item $\mu_\infty = \sup\{ \mu_n: n\in \omega \}$, 
        $\nu_\infty = \sup\{ \nu_n: n\in \omega \}$.   (Note that 
        trivially $\mu_n \le \mu_{n+1} $ and $\nu_n \le \nu_{n+1}$ for
        all $n\in \omega$.)
\end{enumerate}
\end{Definition}

\begin{Fact} \label{easyfact}
Let $L$ be infinite. 
Then $\mu_n(L) \le |L|^+$ and $\nu(L) \le (2^{|L|})^+$. 
\end{Fact}

\begin{Theorem}[Ramsey]\label{ramsey}
For any natural number $k$, $ {\aleph_0}  \to ({\aleph_0} )^2_k$.  
% That is, if 
% $f$ is a function with domain $[A ]^2$ ($A$ infinite) and with range of size
% $k$, then there is an infinite 
% set $A' \subseteq A$ such that $f \on [A']^2$ is constant. 
\end{Theorem}
\begin{proof} See  \cite[10.2]{EHMR}.
 \end{proof}

\begin{Theorem}[Erd\H os+Rado]\label{er}
For any infinite  $\kappa$, 
$(2^\kappa )^+ \to (\kappa ^+)^2_\kappa$. 
%   That is, if $f$ is a
% function with domain of size $[A]^2$ where $A$ is of size $>2^\kappa$, 
% then there is a subset $A' \subseteq A$ of size $>\kappa$ such that
% $f$ is constant on $A'$. 
\end{Theorem}
\begin{proof} See \cite[17.11(32)]{EHMR}.
\end{proof}

\begin{Theorem}[Erd\H os+Rado]\label{er2}
\begin{enumerate}
\item 
If $ \kappa$ is an infinite cardinal, $k$ finite, then 
$(2^{<\kappa})^+ \to (\kappa )^2_k$.
\item If $ \kappa $ is a strong limit cardinal, then 
$$ \kappa^+ \to (\kappa)^2_4$$
\end{enumerate}
\end{Theorem}

\begin{proof} \cite[15.2]{EHMR} proves a theorem that is stronger than
  (1). 
   (2) is a special case of (1). 
\end{proof}

\begin{Corollary}\label{er-cor}
Let $(L, \le)$ be a  partial order.
\begin{itemize}
\itm a If $\kappa $ is an infinite cardinal, $|L| > 2^\kappa$, then
     $\mu(L)> \kappa$.  (In fact, $\mu(L)>\kappa^+$.)
\itm a' $|L| \le 2^{\mu(L)}$. 
\itm b If $L$ is infinite, then $\mu(L) > {\aleph_0} $.
\itm c If $ \kappa $ is a strong limit cardinal, then 
$ \kappa \le  \mu(L)$ iff $ \kappa \le |L|$. 
\itm d If $ \kappa $ is a strong limit cardinal, then
$ | L | > \kappa $ implies $\mu(L) > \kappa $,
\itm e If $ \kappa $ is a strong limit cardinal, then 
 $\mu(L) = \kappa $ implies $|L|= \kappa $. 
\end{itemize}
\end{Corollary}
\begin{proof}
(a) Write $\rho$ for $(2^\kappa)^+$.
Let $(a_i:i<\rho)$ be distinct elements of $L$, and define
$F:[\rho]^2 \to \RR$ by requiring $F(i,j)=R(a_i, a_j)$ whenever $i<j$.
The Erd\H os-Rado theorem \ref{er} promises us an $F$-homogeneous set
$\{i_{\zeta}: {\zeta} < \kappa^+\}$  of size
$\kappa^+ $, which will naturally induce a uniform set
$\{a_{i_{\zeta}}: {\zeta} < \kappa^+\}$ of the same cardinality. 

(a') follows from (a). 

The proofs of (b) and (d) are similar, using \ref{ramsey} and
\ref{er2}, respectively,   instead of \ref{er}. 

(c) follows easily from \ref{easyfact}.    (e) follows from (c) and (d).

\end{proof}

\begin{Canonization Theorem}[Erd\H os+Hajnal+Rado]\label{canon}
Let $\lambda $ be an infinite cardinal,  $(A_i : i < \lambda)$ be a
family of pairwise disjoint sets.
Let $(\kappa_i: i < \lambda)$ be infinite cardinals satisfying
$2^{\kappa_i} < 2^{\kappa _j}$ whenever $i < j$, and assume $|A_i| >
2^{\kappa_i} $.  Let $f$ be a function with domain $[A]^2 $, where $A
= \bigcup _{i< \lambda } A_i$, and let the range of $f$ be small (say:
finite).
Moreover, assume $\kappa_0 \ge 2^{\lambda }$. 

For $\alpha \in A$ write $i_\alpha$ for the unique $i$ such that
$\alpha \in A_i$. 

Then there are sets $(A_i': i < \lambda)$, $|A'_i|= \kappa_i$,
$A'_i \subseteq A_i$,  such
that for $\alpha, {\beta}\in A' $ ($:= \bigcup_{i<\lambda} A'_i$),
$f(\alpha, {\beta})$ depends only on $i_\alpha$ and $i_{\beta}$.  

That is, there is a function $F$ with domain $ [\lambda]^2\cup \lambda $
 such that
for all $\alpha\not= {\beta}$ in $A'$,
 $f(\alpha, {\beta}) = F(i_\alpha,i_{\beta})$.
\end{Canonization Theorem}

\begin{proof} See \cite[28.1]{EHMR}.
\end{proof}

% 
% \begin{Corollary}\label{canon2}
% Assume that $(A_i: i< \lambda)$, $(\kappa_i: i < \lambda)$ and $f$ are
% are as above, and assume moreover that $\lambda = \omega$.    
% 
% Then there is an infinite subset $I \subseteq \omega $, and a constant
% $c$ such that whenever $n\not= m$ are in $I$, $\alpha \in A_n$,
% ${\beta} \in A_m$, then $f(\alpha, {\beta}) = c$.  
% \end{Corollary}
% 
% \begin{proof} Obtain $F$ as in the canonization lemma \ref{canon}, and
% apply Ramsey's theorem \ref{ramsey} to it. 
%\end{proof}

\begin{Corollary}\label{canon2}
  Let $\lambda = cf(\kappa) < \kappa$, $\kappa$ a
strong limit cardinal, $|L| = \kappa$, $(L, \le)$ a 
partial order. 
Then there is is a family 
$(A_{\zeta} : {\zeta} < {\lambda} )$ of subsets of $ L$ 
satisfying 
\begin{enumerate}
\item $(A_\zeta: \zeta < \lambda) $ is canonical. (See
  \ref{defis}(4).)
\item The sequence $(|A_{\zeta} |: {\zeta} < {\lambda})$ is strictly
  increasing.
\item $\sup\{ |A_{\zeta}| : {\zeta} < {\lambda} \}  = \kappa $. 
\item For all $ {\zeta} < {\lambda}$: $|A_{\zeta} | > {\lambda} $. 
\end{enumerate}
\end{Corollary}
\begin{proof}
Find an increasing  sequence 
$(\kappa_i:i<\lambda)$ of cardinals such that 
$2^{\kappa_i}< \kappa_{i+1}$ for all $i$. 
Let $ 
(a_\alpha: \alpha < \kappa)$ be distinct elements, and 
let $A_i:= \{ a_\alpha: \kappa_{i+3} < \alpha < \kappa_{i+4}\}$ for 
$i < \lambda$. 
Thus, $A:= \bigcup_{i < \lambda} A_i$ is a disjoint union, 
and  $|A_i| = \kappa_{i+4}>2^{\kappa_i}$.   
(Note that in this enumeration each set $A_i$ comes 
``before'' $A_j$ for $i<j$.)

Define $f: [A]^2 \to \RR$ by letting
 $f(\{a_\alpha, a_\beta\}) = R(a_\alpha, a_\beta)$
for $\alpha < \beta $, and apply the canonization theorem \ref{canon}.

The resulting sets $(A'_i:i< \lambda) $ will be canonical. 
\end{proof}

\begin{Remark}
If $(\kappa_i: i < \lambda)$ is increasing with limit $\kappa$, then 
$\prod_i \kappa_i = \kappa^\lambda$.   If moreover (as in our case)
$\kappa$ is a strong limit cardinal, then $\kappa^\lambda = 2^\kappa$.
\end{Remark}
\begin{proof} See \cite[6.4]{J}
\end{proof}

% 
% 
% \begin{Definition} For infinite cardinals $\kappa$, $\lambda$ let
% $$L_\lambda(\kappa )= \min \{ \rho:\exists \lambda_0 <\lambda \,\,
% \lambda_0^\rho \ge \kappa \}.$$ 
% \end{Definition}
% 
% 
% 
% \begin{Theorem} \label{inacc}
% \begin{enumerate}
% \item $\kappa^+ \to (\min(\lambda, L_\lambda(\kappa)))^2_\omega $.
% \item If $\kappa$ is a strong limit cardinal, then 
% $\kappa^+ \to (\kappa)^2_\omega$. 
% \end{enumerate}
% \end{Theorem}
% 
% \begin{proof}
% For (1), see \cite[??]{EHMR}.  For (2), note that $L_\kappa(\kappa^+)
% = \kappa$. 
% \end{proof}
% 

\section{Partial maps}

We want to give lower estimates for $\nu(L)$, and then translate
them to lower estimates for $\mu(L)$. (See \ref{munudef} for the 
definitions of $\mu$ and $\nu$.)

Since we will typically only construct many {\em partial}
functions that are pairwise incomparable, we have to give a
sufficient condition that allows us to extend partial
monotone functions to total monotone functions.

\begin{Fact}\label{extend}
Let $L_1$, $L_2$ be partial orders.
 If $f: L_1 \to L_2$ is a partial monotone function whose
range is
%either finite or
 contained in a   complete partial order
$L_2' \subseteq L_2$,
then $f$ can be extended to a total monotone
function $\hat f: L_1 \to L_2$. 
\end{Fact}
\begin{proof}
Let $\hat f(x) = \sup_{L_2'} \{ f(y): y \in \dom(f), y \le x\}$.
\end{proof}

\begin{Corollary} \label{abar}
Let $L$ be a partial order, $\bar A \subseteq L$ a complete partial
order.  (Note: We only require that least upper bounds exist in $(\bar
A, \le)$, we do not care if these bounds are also least upper bounds
in $L$.)

Then for any $A \subseteq \bar A$ we have $\nu(L) \ge \nu(A)$. 
\end{Corollary}

\begin{proof}
Every monotone map $f:A \to A$ can be extended to a monotone map $\hat
f: L \to \bar A$. If $f,g$ are incomparable, then so are  $\hat f$,
$\hat g$. 
\end{proof}

So we will show that $\nu(L)$ is large by showing that 
$\nu(A)$ is large, for some sufficiently ``nice'' $A$.

In our treatment, ``nice'' means in particular ``complete''
(as a partial order), or at least ``contained in a complete p.o.''
Here the following lemma, due to Kaiser and Sauer \cite{KS}
 will be helpful:

\begin{Lemma}
If $(L, \le)$ is an o.p.c.\  lattice, then $L$ is bounded (i.e., has a
greatest and a smallest element). 
\end{Lemma}
\begin{proof} See \cite{KS}.
\end{proof}

Our method to make $\nu(A)$ large will be multiplication:   If
$A_1$, $A_2$, \dots\ are sufficiently ``independent'' (in a
sense to be made precise below), and $f_i:A_i \to A_i$ are
monotone, then we will show that they can be combined to give
a monotone function from $A:=\bigcup_i A_i$ to $A$.

\begin{Independence Lemma}
\label{independent}
Let $L$ be a partial order,
 $A = \bigcup_{i <\lambda } A_i \subseteq L$
and assume that
$(A_i:i < \lambda )$ is canonical. Then:
\begin{enumerate}
\item Whenever $(f_i:i <\lambda )$ is a family of functions, each 
 $f_i:A_i \to A_i$  monotone, then the function 
 $f:= \bigcup_{ i <\lambda } f_i$ is a  monotone function from $A$ to
$A$. 
\item If $B$ is a partial order,
 $(f_i:i <\lambda )$ is a family of functions, each
 $f_i:A_i \to B$  monotone, and if $\bigcup_i A_i$
 is an antichain, then $\bigcup_i f_i$ is monotone from $A$ to $B$.
\item Moreover, if $(f_i:i<\lambda)$ and $(f'_i:i<\lambda)$ are both as
in (1) or (2), and for some $j$ the functions $f_j$ and $f'_j$ are
incomparable, then also $\bigcup_i f_i$ and $\bigcup_i f'_i$ are
incomparable.
\item If $\kappa_i< \nu(A_i)$ for $i<\lambda$, then 
   $ \nu(A) > \prod_{ i <\lambda } \kappa_i$.
\item If $\kappa_i< \nu(A_i,B)$ for $i<\lambda$, and
   $A = \bigcup_{i<\lambda} A_i$ is an antichain, then
   $ \nu(A,B) > \prod_{ i <\lambda } \kappa_i$ (where $B$ is an
   arbitrary partial order)
\end{enumerate}
\end{Independence Lemma}

\begin{proof}
(1) Let $F:I \times I \to \RR$ witness that 
$(A_i:i\in I)$ is canonical. To check that $f$ is monotone, consider an
arbitrary pair $a\le b$ in $A$. 

Now either
there is a single $i$ with $a,b \in A_i$ then $f(a)\le f(b)$
(because $f\on A_i = f_i$ is monotone), or we have $i\not=j$, $a\in A_i$,
$b\in A_j$.   But then we must have $F(i,j)={\boldsymbol<}$, so 
(since $f(a)\in A_i$, $f(b)\in A_j$), we again have 
$f(a) \le f(b)$. 

(2) and (3) are easy. 

(4) follows from (1) and (3), and (5) follows from (2) and (3).
\end{proof}

Now that we know how to get pairwise incomparable functions 
by multiplication, we have to look more closely at the factors
in this product.  The factors are of the form $\nu(A)$, where
$A$ is a uniform set.    
The computation of this cardinal characteristic  turns
out to be easy:

\begin{Fact}\label{..}
\begin{enumerate}
\itm a If $A$ uniform, $|A| > 2$, then
$\nu(A) >2$.
\itm b If $A$ is uniform, $|A| = \kappa \ge  {\aleph_0}$,
  then $\nu(A) > 2^\kappa $, i.e., there are   $2^\kappa$ many pairwise
  incomparable monotone functions from $A$ to $A$.
\itm c If $A$ is an antichain, $|A| = \kappa\ge \aleph_0$, then 
$2^\kappa < \nu(A, \{0,1\})$, i.e.,  there
are $2^\kappa$ many incomparable (necessarily monotone) functions
from $A$ into the two-element lattice $\{0,1\}$.
\end{enumerate}
\end{Fact}
\begin{Proof}
(a) Left to the reader.

(b) This is is trivial if $A$ is an antichain.  
So wlog    assume that $A$ is  well-ordered. 
 Write $A$ as a union of $\kappa$ many 
disjoint  convex sets  $\bigcup_{i<\kappa} A_i$,
each $A_i$ of cardinality $>2$.  Then $(A_i:i \in I)$ is
canonical. So we can apply (a) and 
the independence lemma \ref{independent}(4), and 
get $\nu(A) > 
\prod_{i<\kappa } 2 = 2^\kappa $.

(c): Let $A = \bigcup_i A_i$, where each $A_i$ is of size $> 2$, and
the $A_i$ are pairwise disjoint.  Clearly $\nu(A_i, \{0,1\}) > 2$,
so by the independence lemma \ref{independent} $\nu(A, \{0,1\}) > 2^\kappa$.

\end{Proof}

\section{$\mu$ and $\nu$}

We now turn our attention to the number $\mu_n(L)$.

\begin{Fact} \label{dot}
If $A \subseteq L^n$ is well-ordered of order type $ \kappa $, then
there is $A' \subseteq L$, also well-ordered of order type $ \kappa $.
\end{Fact}
\begin{proof} Let $\bar a^i=(a^i(1), \ldots, a^i(n))$ for $i< \kappa $,
and $i<j  \Rightarrow \bar a^i < \bar a^j$.
For each $k\in \{1, \ldots, n\}$
the sequence $(a^i(k): i< \kappa)$ is weakly increasing.  If the
sequence  $(a^i(k): i< \kappa)$ does not contain a strictly increasing
sequence of length $\kappa$, then it must be eventually constant.
However, this cannot happen for every  $k\in \{1,\ldots, n\}$.
\end{proof}

%That was only the easy case.   If the largeness of
%$ \mu_n(L)$ is
%witnessed by an antichain rather than by a chain, we have to
%use stronger weapons to show that also $\mu(L)$ will be
%large.
%
%\begin{Lemma} \label{mumu}
%Let $(L,\le)$ be a partial order, $\kappa$ an infinite cardinal.
%If $2^{2^\kappa} < \mu_m(L)$, then $\kappa < \mu(L)$.
%\end{Lemma}
%\begin{proof}
%If $(\bar a_i:i <2^{2^\kappa })$,
%$\bar a_i = (a_i(1), \ldots, a_i(n))$,  is a chain in $L^n$, then
%we use \ref{dot} to get a chain of the same length in $L$.   So assume
%that
% $(\bar a_i:i <2^{2^\kappa })$ is an antichain.
%
%For $i<j<2^{2^\kappa}$ let $F(i,j) =
%\bigl( R(a_i(1), a_j(1)), \ldots, R(a_i(n), a_j(n))\bigr) $. Using the
%Erd\H os Rado theorem \ref{er} we get a set $A \subseteq
%2^{2^\kappa}$
%of size $\kappa$ such that $F\on [A]^2$ is constant, so wlog for all
%$i<j<\kappa $ have we
%have $$\bigl( R(a_i(1), a_j(1)), \ldots, R(a_i(n), a_j(n))\bigr) =
%(r_1, \ldots, r_n).$$   Pick $k^* \in \{1, \ldots, n \}$ such that
%$r_{k^*} \not=  {\boldsymbol =}$.  Then $\{ a_i(k^*): i <\kappa \}$ is
%a uniform set of size $\kappa $.
%\end{proof}

Now we finally investigate the relation between 
$\mu$ and $\nu$.  It turns out to be slightly simpler if we
look at $\mu_\infty $ and $\nu_\infty$ first.
%  (Later we 
% will prove $\mu=\mu_\infty$, $\nu=\nu_\infty$, at least for o.p.c.\
% lattices.) 

First we show in \ref{numu} that the existence of many incomparable
monotone functions from $L^n$ to $L$ ($\kappa < \nu_n(L)$) implies the
existence of a large antichain in some $L^m$ $(\kappa < \mu_m(L)$),
assuming that $L$ is o.p.c.  (This is actually the only place in the
whole proof where we talk about lattices rather than general partial
orders.)

Then we show in \ref{munu} that a large (anti)chain in  $L^m$
($\kappa < \mu_m(L)$) implies the existence of $*$very$*$ many
incomparable monotone functions from $L^m$ to $L$ ($2^\kappa < \nu_m(L)$).

These two lemmata are (with minor modifications) taken from 
\cite{HP97}.

Finally in \ref{conclusion} we  combine \ref{numu} and \ref{munu}
to show that $\mu = \mu_\infty $ must be a strong limit cardinal.

\begin{Lemma} \label{numu}
Let $(L,\le)$ be an o.p.c.\  lattice, $\kappa$ a cardinal of uncountable
        cofinality.    If $\kappa < \nu_n(L)$, then $\kappa <
        \mu_\infty(L)$. 
   In particular:
\begin{itemize}
\itm A  $\nu_\infty\le \mu_\infty^+$
\itm B  $\nu_\infty\le \mu_\infty$, if $\mu_\infty $ has uncountable
        cofinality.
\end{itemize}
\end{Lemma}

\begin{proof}
Assume $\kappa <\nu_n(L)$. 
Let $(f_i: i< \kappa)$ be pairwise incomparable functions from $L^n$
to $L$.    Since $L$ is o.p.c., each of these functions is  a polynomial
function.  Thus, for each $i$ there is some natural number $k_i$ and a
lattice-theoretic term $t_i(x_1, \ldots, x_n, y_1, \ldots, y_{k_i})$
and a $k_i$-tuple $\bar b^i=(b^i_1, \ldots, b^i_{k_i})$ such that for all
$a_1$, \dots, $a_n$ we have $f_i(a_1, \ldots, a_n) = 
t_i(a_1, \ldots, a_n, b_1, \ldots, b_{k_i})$. 

Since there are only countably many pairs $(t_i, k_i)$ and we have
assumed $cf(\kappa  )> {\aleph_0}$, we may assume
that they all are equal, say to $(t^*, k^*)$.
  But then $(\bar b^i:i<\kappa)$ must be pairwise
incomparable in $L^{k^*}$, because $\bar b^i \le \bar b^j$ would imply
$f_i\le f_j$. 
Hence we have found an antichain of size $\kappa$ in
$L^{k^*}$.

To get (A), let $\kappa := \mu_\infty^+$, so $cf(\kappa)>\aleph_0$ and
therefore ``$\kappa<\nu_\infty$'' is impossible.  To get (B), let
$\kappa=\mu_\infty$.
 
\end{proof}

From now on we can forget about lattices as long as we only consider
partial orders having properties \ref{numu}(A) and \ref{numu}(B).

\begin{Lemma} \label{munu}
Let $(L,\le , 0, 1)$ be a bounded partial order, $\kappa $ an
infinite cardinal.  If 
$\kappa < \mu_n(L)$, then $2^\kappa < \nu_n(L)$.   In particular, 
$\kappa <\mu_\infty $ implies $2^\kappa <\nu_\infty$. 
\end{Lemma}
\begin{proof}
Let $A \subseteq L^n$ be uniform of size $\kappa$.

Case 1:  $A$ is a chain, so by \ref{dot} wlog $n=1$. Let $\bar A= A \cup
\{0,1\}$. By fact \ref{..}, $\nu(A) >2^\kappa$.   Since $\bar
A$ is a complete partial order,  we may apply fact 
\ref{extend} to get $\nu(L) \ge \nu(A)$. 
Hence $\nu(L)>2^\kappa $. 

Case 2:  $A$ is an antichain.  Use \ref{..}(c). 

%Let $A \subseteq L^n$ be uniform of size $\kappa$, $\bar A = A \cup
% \{0,1\}$.   Clearly there are
% $2^\kappa $ many pairwise incomparable monotone functions [check!!]
% from $ A $ to $\bar A$ , and $\bar A$ is a complete partial order, so
% there are $2^\kappa $ many pairwise incomparable monotone functions
% from $A$ to $L$. 
\end{proof}

\begin{Conclusion}\label{conclusion}    
  If $L$ is infinite and o.p.c., (or, slightly more generally,
if $L$ is an infinite bounded partial order satisfying the conclusion
of \ref{numu}), then 
\begin{itemize}
\itm a $\mu_\infty(L)$ must be a strong limit
        cardinal, 
\itm b $\mu(L)=\mu_\infty(L)$
\itm c $|L| = \mu(L)$. 
\itm d $\nu(L)=|L|$. 
\end{itemize}
\end{Conclusion}

\begin{proof}  
\begin{itemize}
\itm a
 If $\kappa < \mu_\infty(L)$, then $2^\kappa < \nu_\infty(L)$ by
\ref{munu}.   Now  $2^\kappa$ always has uncountable cofinality, so 
we get $2^\kappa  < \mu_\infty(L)$ by \ref{numu}. 

\itm b
Assume that $\mu(L)< \mu_\infty(L)$.
   Let $\lambda = 2^{2^{\mu(L)}} < \mu_\infty(L)$. By \ref{er-cor},
   $|L| \le 2^{\mu(L)} < \lambda$, so $\mu_n(L)\le |L|^+
   \le \lambda$ for all
   $n\in \omega$, hence $\mu_\infty(L) \le \lambda$, a contradiction. 
\itm c   Use \ref{er-cor}(e). 
% Let $ \kappa = \mu(L)$.  
% \relax From   
% \ref{er-cor}(c) we conclude $ \kappa \le |L|$, and from 
% \ref{er-cor}(d) we conclude $ \kappa \ge |L|$. 
\itm d \ref{numu}(B) implies $\nu(L) \le  \nu_\infty(L)\le \mu(L)$, 
 and \ref{munu} implies $\mu_1(L) \le \nu_1(L)$. 
\end{itemize}
\end{proof}

\section{The main lemma}

We have already shown that for an o.p.c.\  lattice $L$ the cardinal 
characteristic 
$\mu(L)$ must be a strong limit cardinal.   In this section we show that
$\mu(L)$ must be regular. 

Letting $ \kappa := \mu(L)$ we first show
 that the singularity of  $ \kappa $ would imply the 
existence of $\gg \kappa $ many incomparable monotone 
functions, and then show
that this would imply $\mu(L)>\kappa $. 

\begin{Main Lemma}\label{mainlemma}
 Let $(L, \le,0,1)$ be a bounded 
partial order, and let $\kappa $ be a
singular strong limit cardinal, $\kappa \le |L|$. 

Then $\nu(L)>\kappa$. 

If moreover  $cf(\kappa) = {\aleph_0}$, then we get even $\nu(L) >
2^\kappa$. 
\end{Main Lemma}

\begin{proof}

Let $ {\lambda} = cf(\kappa)$. 
The first step in the proof of \ref{mainlemma} is to find a canonical
family $(A_i: i < {\lambda})$ which is large, i.e., 
$|\bigcup_i A_i |= \kappa $,  and  $\prod_{i<{\lambda} } |A_i| =
\kappa^{\lambda}  = 2^\kappa $.   If the set $A:= \bigcup_i A_i$
happens to be a chain or antichain , we easily get $2^\kappa $ many
pairwise incomparable monotone functions.   Using the independence
lemma we will show that already  the canonicity of $A$ is sufficient to
get many monotone functions. 

In step 3 we will exhibit many pairwise incomparable (partial) 
functions from $A $ to $A$, so in step 2 we may have to 
 massage $(A_i:i<{\lambda})$ a bit to guarantee that these functions
 can be extended to total functions on $L$. 

What actually happens in steps 2 and 3 depends on whether $ {\lambda}
$ is countable or not.

\subsection*{Step 1}
Let $A = \{a_i: i < \kappa\}$ be distinct elements of our partial order.
For $i<j<\kappa$ let $f(i,j) = R(a_i, a_j)$. (Recall
\ref{rdef}.)

  By the canonization
theorem (or rather, by its corollary \ref{canon2}) we may (after
thinning out our set $A$) wlog assume
that $A = \bigcup_{\zeta < \lambda} A_\zeta$,
where 
the cardinalities $|A_\zeta|$ are increasing with supremum 
$\kappa$, $\lambda < |A_\zeta| <   \kappa$, and 
$( A_\zeta: \zeta < \lambda) $ is canonical.   
Let $\xi:\kappa \to \lambda$ be such that for all $i < \kappa$,
$a_i \in A_{\xi_i}$. 
So 
there is a  function 
$F:[ \lambda]^{\le 2} \to
\{{\boldsymbol<},{\boldsymbol>},{\boldsymbol\|}\}$
 such that $R(a_i, a_j)
= F(\{\xi_i,\xi_j\})$  for all $i<j < \kappa$. 

We may assume $\{0,1\} \cap A = \emptyset$.  Let $\bar A = A \cup
\{0,1\}$.
% , $A_\zeta = \{a_i:\xi_i = \zeta\}$.   

% We may also assume that $|A_{\zeta} | > \lambda$ for all $ {\zeta} $. 

Note that 
$$\prod_{\xi <\lambda }|A_\xi|= \kappa ^{cf(\kappa)}= 2^\kappa.$$

\subsection*{{Step 2}, case a.}
Let us asume   $\lambda = {\aleph_0} $ for the
moment. So we have a canonical sequence 
$(A_n:n\in \omega)$, witnessed by $F:[\omega]^2 \to \{
\boldsymbol{<},
\boldsymbol{>},
\boldsymbol{\|}\}$. By Ramsey's theorem \ref{ramsey} there is an
  infinite  set   $X \subseteq  \omega $ such that $F$ is constant on
  $[X]^2$. By dropping some elements of the sequence $(A_n:n\in
  \omega)$  [i.e., replacing 
  $(A_n:n\in  \omega)$ by $(A_n:n \in X)$, and then for notational
  simplicity only pretending that $X=\omega$] 
  we may assume that $F$ is  constant, say
$$ \forall n\, \forall k : \ n < k \  \Rightarrow \ 
        \forall a \in A_n \, \forall b\in A_k: 
R(a,b)=c$$

\begin{figure}
\setlength{\unitlength}{0.0125in}%
\begin{picture}(390,210)(75,575)
\thicklines
\put( 85,590){\line( 1, 0){ 30}}
\put( 75,640){\line( 1, 0){ 50}}
\put( 65,690){\line( 1, 0){ 70}}
\put(285,580){\vector( 0, 1){ 25}}
\put(285,630){\vector( 0, 1){ 35}}
\put(285,690){\vector( 0, 1){ 45}}
\put(425,605){\vector( 0,-1){ 25}}
\put(425,665){\vector( 0,-1){ 35}}
\put(425,740){\vector( 0,-1){ 45}}
\put(150,585){\makebox(0,0)[lb]{\raisebox{0pt}[0pt][0pt]{$A_0$}}}
\put(150,635){\makebox(0,0)[lb]{\raisebox{0pt}[0pt][0pt]{$A_1$}}}
\put(150,690){\makebox(0,0)[lb]{\raisebox{0pt}[0pt][0pt]{$A_2$}}}
\put(105,735){\makebox(0,0)[lb]{\raisebox{0pt}[0pt][0pt]{$\vdots$}}}
\put(285,750){\makebox(0,0)[lb]{\raisebox{0pt}[0pt][0pt]{$\vdots$}}}
\put(425,750){\makebox(0,0)[lb]{\raisebox{0pt}[0pt][0pt]{$\vdots$}}}
\put(295,590){\makebox(0,0)[lb]{\raisebox{0pt}[0pt][0pt]{$A_0$}}}
\put(295,640){\makebox(0,0)[lb]{\raisebox{0pt}[0pt][0pt]{$A_1$}}}
\put(295,715){\makebox(0,0)[lb]{\raisebox{0pt}[0pt][0pt]{$A_2$}}}
\put(210,655){\makebox(0,0)[lb]{\raisebox{0pt}[0pt][0pt]{ or }}}
\put(350,660){\makebox(0,0)[lb]{\raisebox{0pt}[0pt][0pt]{ or }}}
\put(440,585){\makebox(0,0)[lb]{\raisebox{0pt}[0pt][0pt]{$A_0$}}}
\put(440,645){\makebox(0,0)[lb]{\raisebox{0pt}[0pt][0pt]{$A_1$}}}
\put(440,720){\makebox(0,0)[lb]{\raisebox{0pt}[0pt][0pt]{$A_2$}}}
\end{picture}
% 
% 
% 
% 
% 
% 
% \setlength{\unitlength}{0.0125in}%
% \begin{picture}(390,210)(65,575)
% \thicklines
% \put( 65,590){\line( 1, 0){ 70}}
% \put( 65,640){\line( 1, 0){ 70}}
% \put( 65,690){\line( 1, 0){ 70}}
% \put(285,580){\vector( 0, 1){ 30}}
% \put(285,630){\vector( 0, 1){ 35}}
% \put(285,700){\vector( 0, 1){ 35}}
% \put(430,605){\vector( 0,-1){ 30}}
% \put(430,665){\vector( 0,-1){ 35}}
% \put(430,740){\vector( 0,-1){ 35}}
% \put(150,585){\makebox(0,0)[lb]{\raisebox{0pt}[0pt][0pt]{$A\_0$}}}
% \put(150,635){\makebox(0,0)[lb]{\raisebox{0pt}[0pt][0pt]{$A\_1$}}}
% \put(150,690){\makebox(0,0)[lb]{\raisebox{0pt}[0pt][0pt]{$A\_2$}}}
% \put(105,745){\makebox(0,0)[lb]{\raisebox{0pt}[0pt][0pt]{$\vdots$}}}
% \put(280,770){\makebox(0,0)[lb]{\raisebox{0pt}[0pt][0pt]{$\vdots$}}}
% \put(425,770){\makebox(0,0)[lb]{\raisebox{0pt}[0pt][0pt]{$\vdots$}}}
% \put(305,590){\makebox(0,0)[lb]{\raisebox{0pt}[0pt][0pt]{$A_0$}}}
% \put(295,640){\makebox(0,0)[lb]{\raisebox{0pt}[0pt][0pt]{$A_1$}}}
% \put(295,715){\makebox(0,0)[lb]{\raisebox{0pt}[0pt][0pt]{$A_2$}}}
% \put(210,655){\makebox(0,0)[lb]{\raisebox{0pt}[0pt][0pt]{ or }}}
% \put(360,660){\makebox(0,0)[lb]{\raisebox{0pt}[0pt][0pt]{ or }}}
% \put(455,585){\makebox(0,0)[lb]{\raisebox{0pt}[0pt][0pt]{$A_0$}}}
% \put(450,645){\makebox(0,0)[lb]{\raisebox{0pt}[0pt][0pt]{$A_1$}}}
% \put(450,720){\makebox(0,0)[lb]{\raisebox{0pt}[0pt][0pt]{$A_2$}}}
% \end{picture}
% 
% 
\centerline{\small Figure 1}
\end{figure}

There are (at most) 9 possible
types of our family  $(A_n:n\in \omega)$: Each 
 $A_n$ can be  well-ordered,
co-well-ordered, or an antichain, and there are also $3$ possible
values for $c$.  For example, if 
$c=$``$\boldsymbol{<}$'', then set $A$ has one of the
3 forms given in the figure~1.

However, our construction of a large family of incomparable 
functions will be ``uniform'', i.e.,
be the same in all cases. 
Note that if 
 $c\in \{ \boldsymbol{<},
\boldsymbol{>}\}$ and each $A_n$ is an antichain (leftmost possibility
in the above picture), then $\bigcup_n A_n$ is not a complete partial
order, since every element of $A_{n+1}$ is a minimal upper (or maximal
lower) bound for any nontrivial subset of $A_n$.

For each $k$ we let 
$A'_{2k+1}$ be a singleton subset of $A_{2k+1}$, 
and we let $ A'_{2k} = A_{2k}$.    Now let 
$ B:= \bigcup_n A'_n \cup \{ 0,1\}$.

\subsection*{Step 3, case a.}
We are still assuming $ \lambda = \aleph_0$. 
 It is easy to see (by considering cases --- one of them 
is sketched in figure~2) 
that the set $B$ defined in step 2 (case  a)
is a complete partial order.    We leave the details to the reader.

\begin{figure}
\begin{center}
\setlength{\unitlength}{0.0125in}%
\begin{picture}(90,250)(150,580)
\thicklines
\put(182,628){$\bullet$}
\put(185,705){$\bullet$} %\circle{8}
\put(150,670){\line( 1, 0){ 75}}
\put(150,590){\line( 1, 0){ 75}}
\put(155,745){\line( 1, 0){ 75}}
\put(240,580){\makebox(0,0)[lb]{\raisebox{0pt}[0pt][0pt]{$ A_0$}}}
\put(235,665){\makebox(0,0)[lb]{\raisebox{0pt}[0pt][0pt]{$ A_2$}}}
\put(240,740){\makebox(0,0)[lb]{\raisebox{0pt}[0pt][0pt]{$ A_4$}}}
\put(200,620){\makebox(0,0)[lb]{\raisebox{0pt}[0pt][0pt]{$ A_1$}}}
\put(200,700){\makebox(0,0)[lb]{\raisebox{0pt}[0pt][0pt]{$ A_3$}}}
\put(195,770){\makebox(0,0)[lb]{\raisebox{0pt}[0pt][0pt]{$\vdots $}}}
\end{picture}
\end{center}
\centerline{\small Figure 2}
\end{figure}

Note that we still have $\prod_{n<\omega} |A'_n| = 
 \prod_{n\in\omega} |A'_{2n}| = \kappa ^ {\lambda} = 2^\kappa$, 
since the cardinalities $(|A_{2n}|:n \in \omega)$ are 
also increasing to $\kappa$. 

 By \ref{abar} and the independence lemma \ref{independent}, 
$\nu(L)\ge \nu(A) > 2^\kappa$. 
This concludes the discussion of the case $ {\lambda} = {\aleph_0} $. 

% are 
% First we claim that whenever 
% $0 < \zeta$, and  $f: A_0 \to A_0$ and $g:
% A_\zeta \to A_\zeta$ are both monotone.   Then $f\cup g$ is also a
% monotone partial function from $\bar A $ to $\bar A$.   [Proof:$(A_0,
% A_\xi)$ is canonical, so the independence lemma \ref{independent}
% applies.] 
% 
% Conclusion: There are at least $\kappa $ many pairwise incomparable
% functions from $\bar A $ to $\bar A$. 
% 

\subsection*{Step 2, case b.}  Now we assume $\lambda >{\aleph_0} $.      
Note that now we are only aiming for $ \kappa $ many pairwise
incomparable monotone functions.  
Again let $F: \lambda \times \lambda \to \{ 
{\boldsymbol<} ,
{\boldsymbol>}   ,
{\boldsymbol\|}    \}$ witness canonicity. 

 By
 dropping some of the sets $A_{\zeta} $ we may assume that the
 function $\zeta \mapsto F(0, {\zeta} )$ is constant, say with some
 value $c_0$.  Thinning out three more times we may assume that there
         are  constants $c_0, c_1, c_2, c_3 \in 
         \{ {\boldsymbol{<}}, {\boldsymbol{>}}, 
        {\boldsymbol{\|}}\}$ such that 
        \begin{itemize}
         \itm a $ \forall {\zeta}  > 0$, $F(0, {\zeta} )= c_0$ 
         \itm b $ \forall {\zeta}  > 1$, $F(1, {\zeta} )= c_1$ 
         \itm c $ \forall {\zeta}  > 2$, $F(2, {\zeta} )= c_2$ 
        \itm d $ \forall {\zeta}  > 3$, $F(3, {\zeta} )= c_3$ 
        \end{itemize}

        Choose $i_0 < i_1 \in \{0,1,2,3\}$ such that $c_{i_0} = c_{i_1}$,
        write $A_*$ for $A_{i_0}$,  
        let $a_*$ be any element of $A_{i_1}$, and let $c_*= c_{i_0}$.  

        Depending on the value of $c_* $ we now have one the following
        possibilities:
        \begin{itemize}
        \itm a  either   every element of $A_*$ is incompatible 
                with every element
        of any $A_{\zeta} $, $ {\zeta} >4$, 
        \itm b or we have for all  $ {\zeta} >4 $:
        $$ \forall x \in A_*\,\,  \forall y \in A_{\zeta} :  x < a_* <y$$
        \itm c or the dual of (b) is  true. 
        \end{itemize}

        \subsection*{Step 3, case b.}

Since  all sets  $A_{\zeta}$ are infinite, we have 
$\nu(A_\zeta)> |A_\zeta|$, so we can find
pairwise incomparable monotone functions 
$(f_{{\zeta}, i}:i \in A_\zeta)$ 
from $A_\zeta $  to $A_{\zeta} $. 

Let $f^*$ be the identity function on $\{a^*\}$, and let 
$(f_{*,\zeta}: 4<{\zeta} < {\lambda} )$ be a family of 
pairwise incomparable monotone functions from $A_*$ to $A_*$.
(Recall that all our sets $A_{\zeta} $ had cardinality $ {\lambda} $,
so it is possible to find that many functions.)

Note that for $4 <{\zeta}$ the family $(A_*, \{a_*\}, A_{\zeta})$
is canonical, so
by the independence lemma  we can conclude 
that $f_{*,{\zeta}} \cup f_* \cup  f_{\zeta, i}$ is a
monotone function. Moreover, $\{0,1 , a_*, \}\cup A_*  \cup A_{\zeta}$ is a
complete partial order (again we leave the easy task of checking this
fact to the reader), so by \ref{extend} the function
$f_{*,{\zeta}} \cup f_* \cup f_{\zeta, i}$ can be extended to a total 
monotone function $\hat f_{{\zeta},i}:L\to L$.

Clearly any two functions $\hat f_{{\zeta},i}, \hat f_{{\zeta} ', i'}$
are incomparable:
If ${\zeta} \not= {\zeta} '$ then this is due to the incomparability
of 
$f_{*,{\zeta} }$ and $f_{*,{\zeta}' }$, and for ${\zeta}={\zeta} '$
we use the incomparability of 
$f_{{\zeta} , i }$ and $f_{{\zeta} , i' }$. 

Note that the cardinality of the index set 
$\{ ({\zeta}, i):  4<{\zeta} <\lambda, i \in A_\zeta \}$ is
$ \kappa $.

This concludes the discussion of the case $ {\lambda} >  {\aleph_0} $,
and hence also the proof of the main lemma. 
\end{proof}

\begin{Conclusion}
If $(L, \le)$ is an o.p.c.\  lattice (or, slightly more generally, if
$L$ is a partial order satisfying the conclusion of \ref{numu}), 
then $\mu(L)=|L|$  is  an inaccessible
cardinal. 
\end{Conclusion}

\begin{proof} Let $\kappa = \mu(L)$.   From \ref{conclusion}
 we know that $\kappa$ is a strong limit, and that $|L|= \kappa $.
   Assume that $\kappa$ is  singular.  

First assume that $cf(\kappa)$  is uncountable.  The main lemma tells
us that $\nu(L) > \kappa$, so by \ref{numu} we conclude $\mu_\infty(L)
> \kappa$, a contradiction. 

Now we consider the second case:  $cf(\kappa) = {\aleph_0}$. Here the
main lemma tells us $\nu(L) > 2^\kappa$.  Since $2^\kappa$ has
uncountable cofinality, we can again apply \ref{numu} and again get
$\mu_\infty(L) > 2^\kappa > \kappa$, a contradiction.

% Finally we have to show $|L|= \kappa$.    But from \ref{inacc} we know
% that $\kappa^+\to (\kappa)^2_4$, so if $|L| \ge \kappa^+$ we would
% have a uniform set of size $\kappa$ in $L$, a contradiction. 

\end{proof}

\begin{Remark}
Note that the cardinality of an o.p.c.\ lattice cannot be a weakly
compact cardinal.  
\end{Remark}
% \begin{proof}  If $\kappa:= |L|$ is weakly compact, then 
% $\mu(L)>\kappa$, so $\nu(L)>2^\kappa$. 
% \end{proof}

\bibliographystyle{lit-unsrt}
\bibliography{listb,listx}

\end{document}